\documentclass{primus}

\usepackage{amsmath,amssymb, amsthm}
\usepackage{color}
\usepackage{graphicx, subcaption}
\usepackage{float}
\restylefloat{table}
\usepackage{hyperref}

\newcommand{\edits}[1]{#1}
\newcommand{\newedits}[1]{#1}

\newcommand{\vect}[1]{\mathbf{#1}}

\title{Principal Component Analysis: Resources for an Essential Application of Linear Algebra}

\author{S. Pankavich\\
Department of Applied Mathematics and Statistics\\
Colorado School of Mines\\
Golden, CO 80403\\
pankavic@mines.edu
(303) 273-3584\\
\and
R. Swanson\\
Department of Applied Mathematics and Statistics\\
Colorado School of Mines\\
Golden, CO 80403\\
swanson@mines.edu\\
}


\newcommand{\AmSLaTeX}{$\cal A$\kern-.1667em\lower.5ex\hbox{$\calM$}\kern-.125em $\cal S$-\LaTeX}

\usepackage[framed,numbered,autolinebreaks,useliterate]{mcode}

\newcommand{\bfR}{\mathbb{R}}

\newtheorem*{theorem*}{Theorem}

\begin{document}


\makePtitlepage
\makePtitle

\begin{abstract}
Principal Component Analysis (PCA) is a highly useful topic within an
introductory Linear Algebra course, especially since it can be
used to incorporate a number of applied projects.
This method represents an essential application 
and extension of the Spectral Theorem and is commonly used within a variety of fields,
including statistics, neuroscience, and image compression.
We present a synopsis of PCA and include a number of examples that can be used within upper-level mathematics courses to 
engage undergraduate students while introducing them to one of the most widely-used applications of linear algebra.
\end{abstract}

\listkeywords

\section{INTRODUCTION}
Though it is typically the first proof-based course most students experience, linear algebra is also \edits{an important topic that is used to solve problems in a wide variety of fields}.
This is why, at many universities, the course is a requirement for a large number of degree programs in engineering and the natural sciences.
However, since many colleges, especially those with a technical focus, possess significantly more of these majors than those within mathematics, the majority of enrolled students often find difficulty with the abstract nature of the subject, especially when the question\edits{s} ``How does this apply to my major?'' \edits{and ``What are the real-world applications?''}, can largely go unanswered.  

One valuable application that can be included in a first Linear Algebra course is the method of Principal Component Analysis (PCA). 
PCA is a good choice for an applied example to which linear algebra is crucial because it arises in so many different contexts, as we will demonstrate within subsequent sections.
The method arises in countless \edits{disciplines}, including \edits{but not limited to} statistics \cite{J}, electrical engineering \cite{PCAEA},
genetics \cite{RSA}, neuroscience \cite{Peyrache}, facial recognition \cite{Turk}, control theory \cite{PCAEA}, and mechanical and systems engineering \cite{PCAME}.
\edits{Perhaps the largest obstacle to including PCA within a first semester Linear Algebra course is that it is often performed using the Singular Value Decomposition (SVD), which is a topic usually reserved for either a second course in the subject, a course in computational linear algebra, or a graduate linear algebra course (see \cite{Remski,Kalman} for more information on the SVD).  Of course, many universities do not offer these follow-on courses, and hence students are not exposed to PCA, even though it is such a widely-utilized tool.  However, we note that PCA can be both understood and performed without knowledge of the SVD, using a topic that is usually included within a first-semester course, namely the Spectral Theorem.  We will elaborate further on how this is accomplished later on. }
\edits{In addition,} though PCA is sometimes included as an applied topic within linear algebra texts, most notably \cite{Lay} and \cite{Strang}, it is
often overlooked as an important application within other sources.  
Even when \edits{discussed in these texts}, it is often difficult to describe the utility or impact of the
method because examples within textbooks \edits{are typically} computed by hand, while real-world implementation is always performed using a computer.


To address these difficulties, we have created resources, both theoretical and computational, for instructional use within the classroom of a first course in linear algebra.  These resources require a computer in the classroom, but will not require students to understand numerical aspects of linear algebra.
We also provide the code that generates the data sets, numerics, and images used within our examples.  Instructors who are comfortable with {\it MATLAB} can alter our code to make the presentation more interactive, while those with less computational interest or experience need not do so in order to generate all of the material we present.
While we use {\it MATLAB} to perform computations, the principles and algorithms can be applied using other (perhaps open-source) software, such as {\it Octave}.
For instructors without access to {\it MATLAB}, we also provide a number of programs on the first author's website:
\begin{center} \url{http://inside.mines.edu/~pankavic/PCAcode/Octave}.
 \end{center}
that are similar to those in the appendix but written using {\it Octave}.

\edits{These computational tools arise from a first-semester Linear Algebra course conducted at the Colorado School of Mines,
and have had a strong impact on the interest of non-majors and, subsequently, the increase in mathematics minors at our institution.  Both authors have implemented portions of this material in their classrooms and have received positive feedback from students regarding their interest in PCA applications, especially concerning its use in image compression and statistics. In addition, the inclusion of PCA has benefitted the mathematical breadth and depth of our applied majors and increased their exposure to new and interesting applications of linear algebra.}
Hence, we believe the resources included herein can be used to introduce students to a variety of real-world applications of the subject without the need for a course in scientific computing.

We also mention that the implementation of these materials is fairly robust.  Instead of a Linear Algebra course, they could be included within an introductory or second course in statistics for students who have already taken linear algebra, or within a scientific computing course, assuming students possess the necessary background. 
In general, this article provides a number of teaching resources and examples to engage students
in learning about this essential application of linear algebra.

\section{BACKGROUND} 

We begin by recalling \edits{a crucial} theorem - one of the most important results in all of linear algebra - on which \edits{we base our discussion of PCA}.  \edits{Before stating the result, we note that a stronger version of the theorem extends to a more general class of matrices over the complex numbers. } 

\begin{theorem*} (Spectral Theorem \cite{Lay})
Let $n \in \mathbb{N}$ be given and $A$ be an $n \times n$ matrix of real numbers with $A^T$ denoting its transpose. \edits{Then, $A$ is symmetric (i.e., $A^T = A$) if and only if $A$ is orthogonally diagonalizable.  In this case, $A$ possesses $n$ real eigenvalues, counting multiplicities.} 
\end{theorem*}
This is an extremely powerful result and precisely guarantees \newedits{for such A} the existence of
$\lambda_k \in \bfR$ and orthonormal column vectors $\vect{v}_k \in \bfR^n$ for every $k=1,...,n$ such that
\begin{equation}
\label{ST}
A = \lambda_1 \vect{v}_1  \vect{v}_1^T + \cdots +  \lambda_n \vect{v}_n \vect{v}_n^T
\end{equation}
arising from the orthogonal diagonalization.
Because the $\vect{v}_k$ vectors are orthonormal, each associated matrix $\vect{v}_k\vect{v}_k^T$ \newedits{is orthogonal},
and thus $A$ can be decomposed into a sum wherein each term is an eigenvalue multiplied by a rank one matrix \newedits{generated by a unit vector}.
Hence, the eigenvalues alone determine the magnitude of each term in the sum, while the eigenvectors
determine the directions. 
These eigenvectors are called \emph{principal components} or \emph{principal directions} and we will expand upon this further in the next section.
We begin the discussion of the specifics of PCA with an introductory example.

\subsection{Introductory Height \& Weight Problem}

Consider a study in which we want to determine whether or not the heights and weights of a group of individuals are correlated.  That is, we want to know
whether the known value of a person's height seems to dictate whether they tend to be heavier or lighter, and thus influences their weight.  Assume we are given data for $30$ specific people, displayed within Table~\ref{tab:1}.  For this example, our data set originates from a commonly available study \cite{SOCR}.  In practice and for large enough class sizes, 
the data can be obtained in the classroom, directly from students.  For instance, \edits{during} the class period before implementation of this example, the instructor can collect anonymous information regarding student height and weight, compile and store this data, and integrate it into the project for the following class period. \edits{As an alternative, this problem could serve as a project for a group of students within the course, wherein the group collects the data, performs PCA using the code provided in the appendix, interprets the results, and presents their findings in a brief class seminar.} For those with smaller class sizes or concerns regarding anonymity, data can be taken from our source \cite{SOCR} that is readily available on the \newedits{Internet}.
\begin{table}[h]
\begin{tabular}{|l |*{6}{c}|}
\hline
Person & 1 & 2 & 3 & 4 & 5 & 6\\
\hline
\hline
Height & 67.78 & 73.52 & 71.40 & 70.22	& 69.79 & 70.70\\
Weight & 132.99 & 176.49 & 173.03 & 162.34 & 164.30 & 143.30\\
\hline
Person & 7 & 8 & 9 & 10 & 11& 12\\
\hline
\hline
Height & 71.80& 72.01& 69.90 & 68.78 & 68.49 & 69.62 \\
Weight  & 161.49 & 166.46 & 142.37 & 150.67 & 147.45 & 144.14\\
\hline
Person & 13 & 14 & 15 & 16 & 17 & 18\\
\hline
\hline
Height & 70.30 & 69.12 & 70.28 & 73.09 & 68.46 & 70.65\\
Weight  & 155.61 & 142.46 & 146.09 & 175.00 & 149.50 & 162.97\\
\hline
Person & 19 & 20 & 21 & 22 & 23 & 24\\
\hline
\hline
Height & 73.23 & 69.13 & 69.83 & 70.88 & 65.48 & 70.42\\
Weight  & 177.90 & 144.04 & 161.28 & 163.54 & 126.90 & 149.50\\
\hline
Person & 25 & 26 & 27 & 28 & 29 & 30\\
\hline
\hline
Height & 69.63 & 69.21 & 72.84 & 69.49 & 68.53 & 67.44\\
Weight  & 161.85 & 149.72 & 172.42 & 151.55 & 138.33 & 133.89\\
\hline
\end{tabular}
\centering \caption {  \footnotesize Heights (in.) and Weights (lbs.) for $30$ young adults \cite{SOCR}.}
\label{tab:1}
\end{table}

\begin{figure}[t]
\centering
\includegraphics[scale=0.4]{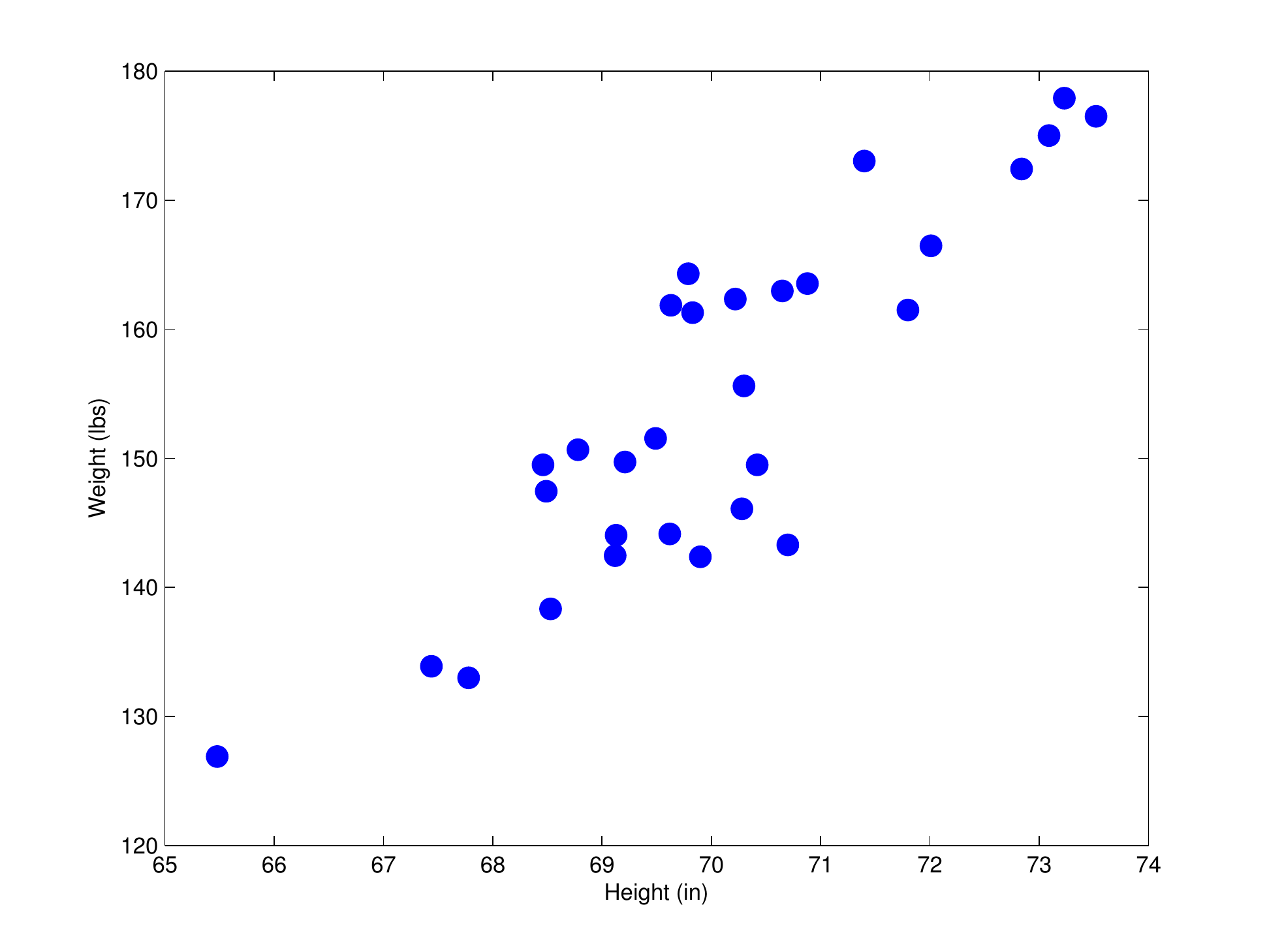}
\vspace{-0.2in}
\caption{ \label{fig:data}  \footnotesize  \edits{Plot} of Height/Weight \edits{datapoints}.}
\end{figure}

\begin{figure}[t]
\begin{subfigure}[t]{0.45\textwidth}
\hspace{-0.3in}
\includegraphics[scale=0.33]{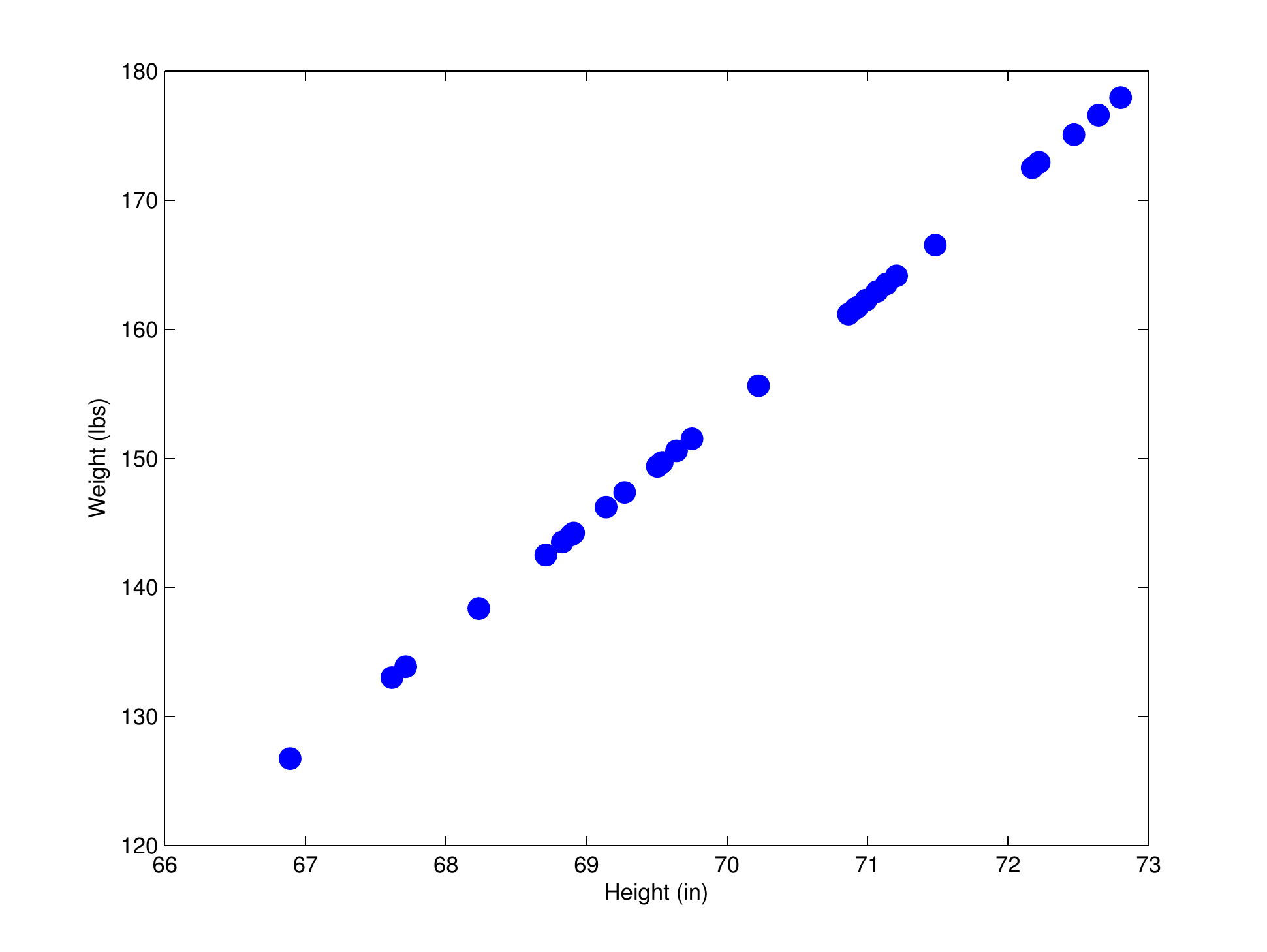}
\vspace{-0.2in}
\end{subfigure}
\hspace{0.1in}
\begin{subfigure}[t]{0.45\textwidth}
\includegraphics[scale=0.33]{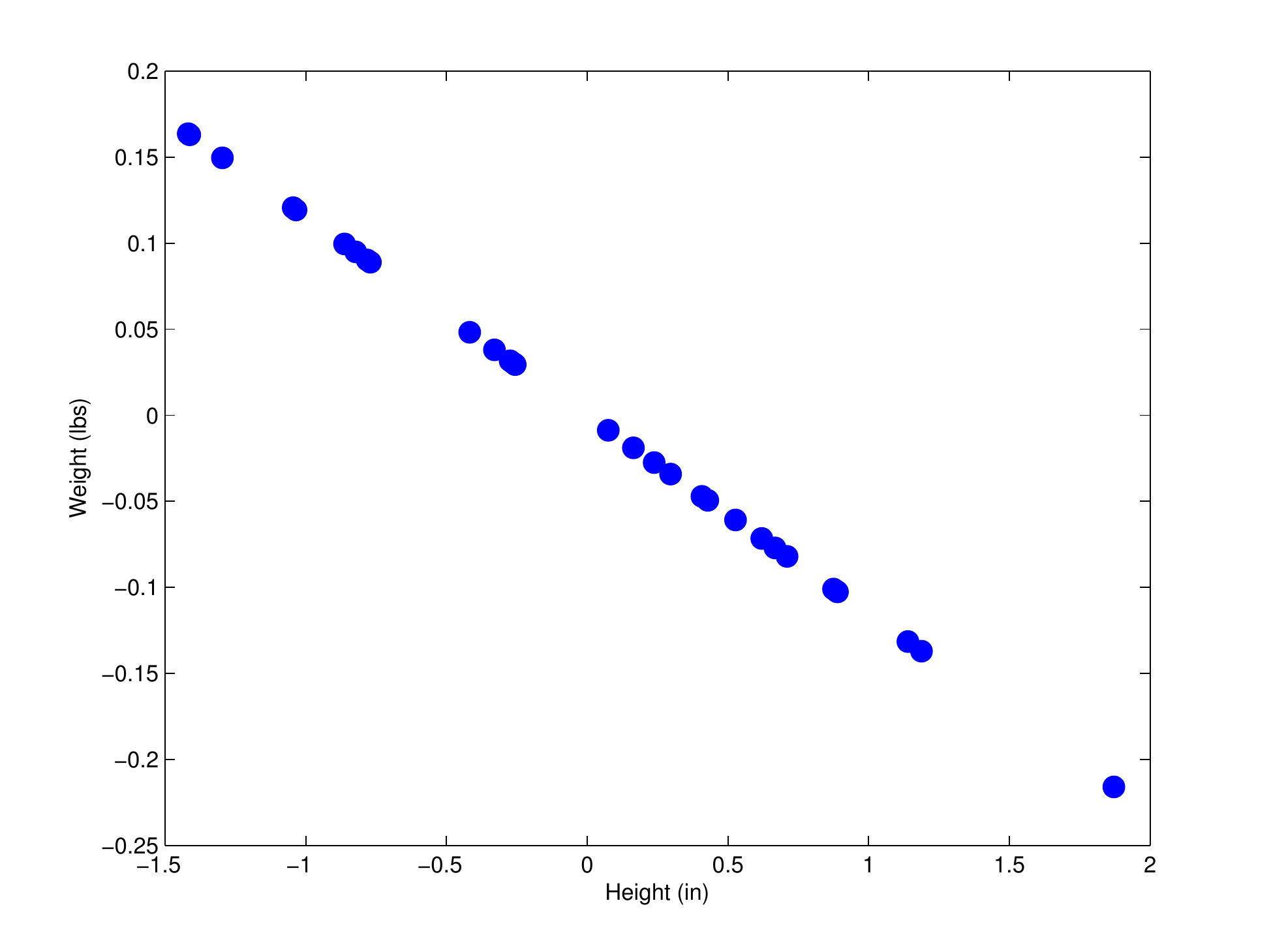}
\vspace{-0.2in}
\end{subfigure}
\caption{ \label{fig:PCintro} \footnotesize  Height/Weight data projected onto the principal components (left - $v_1$; right - $v_2$) .  By the Spectral Theorem, the data represented in Fig.~\ref{fig:data}
is exactly the sum of the projections onto these two components.}
\end{figure}

Since the question of interest is whether the two measured variables, height and weight, seem to change together, the relevant quantity to consider is the covariance of the two characteristics within the data set. 
This can be formed in the following way.  
First, the data is stored in $X$, a $2 \times 30$ matrix. 
\edits{Then, the entries are used to compute the mean in each row, which will be used to center \edits{or ``mean-subtract''} the data}. \edits{This latter step is essential, as many of the results concerning PCA are only valid upon centering the data at the origin.}  Computing the means of our measurements (Table~\ref{tab:1}), we find 
$$\mbox{\boldmath${\mu}$} = \left[\begin{array}{r} 70.06 \\ 154.25 \end{array}\right].$$
Using $x_{ij}$, the entries of the data matrix $X$, the associated $2 \times 2$ covariance matrix $S$ is constructed with entries
$$s_{ik} = \frac{1}{30-1} \sum_{j=1}^{30} (x_{ij} - \mu_i) (x_{kj} - \mu_k)$$
so that
$$S = \left[\begin{array}{rr} 3.26 & 21.72 \\ 21.72 & 188.96 \end{array}\right]. $$
Notice that this matrix is necessarily symmetric, so using the Spectral Theorem it can be orthogonally diagonalized.  Upon computing the eigenvalues and eigenvectors \newedits{of S}, we find $\lambda_1 = 191.46$, $\lambda_2 = 0.76$, and 
$$ \vect{v}_1 = \left[\begin{array}{r} 0.11 \\ 0.99 \end{array}\right]
\qquad 
\vect{v}_2 = \left[\begin{array}{r} -0.99 \\ 0.11 \end{array}\right]. $$
Here, $\vect{v}_1$ and $\vect{v}_2$ are the \emph{principal components} of the \edits{covariance matrix $S$ generated by the data matrix $X$,}
as previously described.
Thus,  we see from (\ref{ST}) that
$$S = \lambda_1 \vect{v}_1 \vect{v}_1^T + \lambda_2 \vect{v}_2 \vect{v}_2^T$$
and because the difference in eigenvalues is so large, it appears
that the first term is responsible for most of the information encapsulated
within $S$. 
Regardless, we can re-express the given data in the new orthonormal basis generated by $\vect{v}_1$ and $\vect{v}_2$
by computing the coordinates $P^T X$ where
$$P = \left[\begin{array}{rr} 0.11 & -0.99 \\ 0.99 & 0.11 \end{array}\right]$$
is the orthogonal matrix whose columns are $\vect{v}_1$ and $\vect{v}_2$.
In fact, we could left multiply the data matrix by each component separately, namely $\vect{v}_1^T X$
and $\vect{v}_2^T X$, to project the data onto each principal direction (Fig.~\ref{fig:PCintro}).
Hence, the data can be separated into projections along $\vect{v}_1$ and $\vect{v}_2$, respectively. 
\edits{We see from looking at the scales} in Figure~\ref{fig:PCintro} that the heights and
weights along $\vect{v}_2$ are significantly less than those along $\vect{v}_1$, which tells us that the majority of the information contained within $X$ lies along $\vect{v}_1$.  
Computing the slope of the line in the direction of $\vect{v}_1$ and choosing a point thru which
it passes, we can represent it by
$$y - 154.25 = 9 (x - 70.06),$$ where $x$ represents the height of a given individual and $y$ is 
their corresponding weight.
Hence, we see that height and weight appear to be strongly correlated, 
and PCA has determined the direction with optimal correlation between the variables.

\begin{figure}[t]
\centering
\includegraphics[scale=0.4]{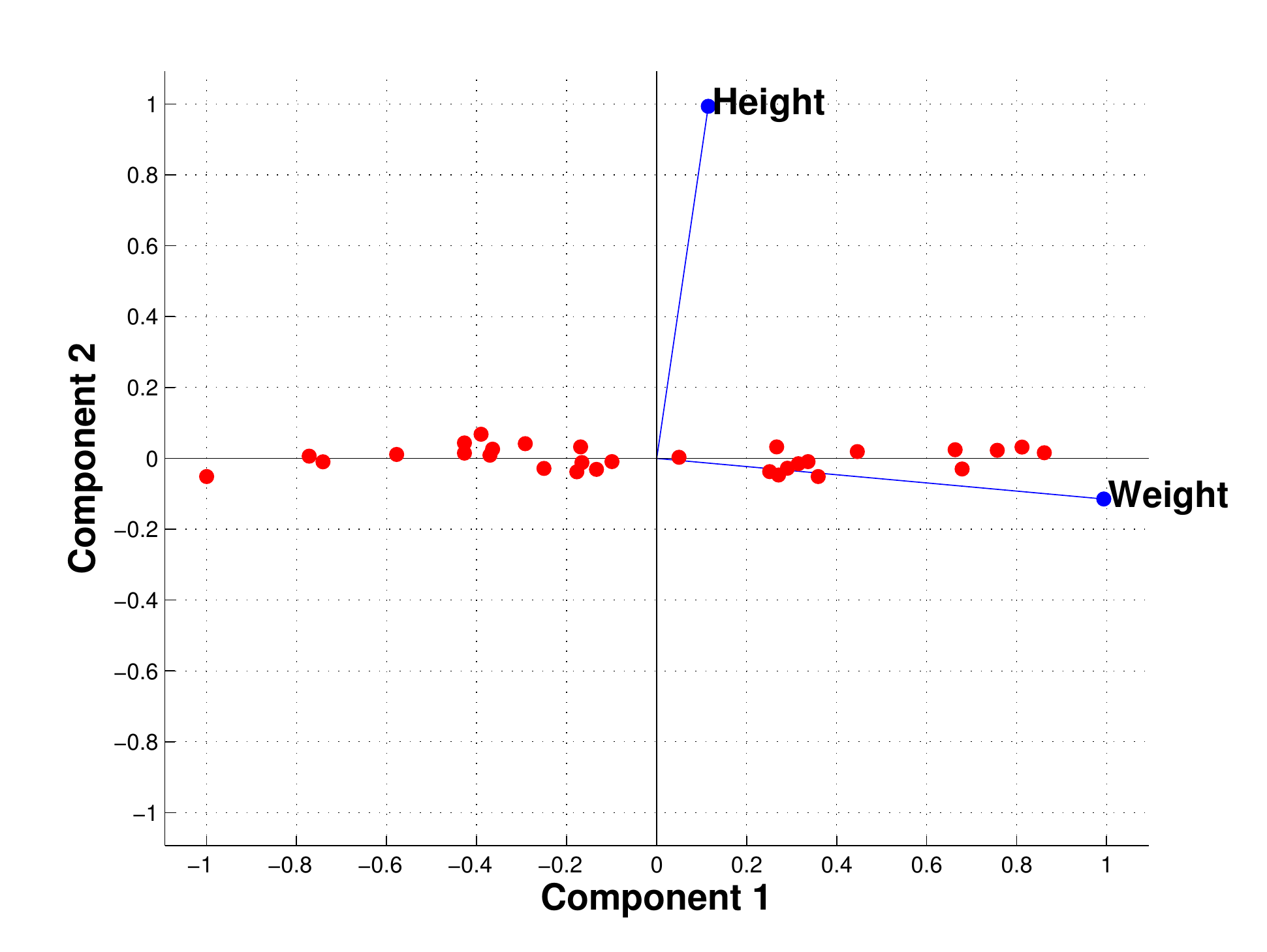}
\vspace{-0.2in}
\caption{ \footnotesize \label{fig:HtWt_biplot} Biplot of Height/Weight data with $2$ principal components.  \edits{The blue Height and Weight vectors are displayed as linear combinations of the principal components.} Note that \edits{these principal components} effectively rotate the \edits{height and weight} data in the plane.}
\end{figure}

The principal component analysis for this example took a small set of data and identified
a new orthonormal basis in which to re-express it. 
In two dimensions the data are effectively rotated to lie along the line of best fit (Fig.~\ref{fig:HtWt_biplot}), with the second principal
direction merely representing the associated unit orthogonal complement of the first.
This mirrors the general aim of PCA: to obtain a new orthonormal basis that organizes the data optimally, in the sense that the variance contained within the vectors is maximized along \edits{successive} principal component(s).

\subsection{Summary of PCA}
\label{sec:steps}

In short, PCA can be performed to compute an optimal\newedits{, ordered} orthonormal basis of a given set of vectors, or data set, in the following steps.
\begin{enumerate}
\item Gather $n$ samples of $m$-dimensional data, i.e. vectors $\vect{d}_1, ..., \vect{d}_n \in \bfR^m$ stored in the
\edits{$m \times n$} matrix $X$ with columns $\vect{d}_1, ..., \vect{d}_n$, so that $x_{ij}$ represents
the $i^{th}$ entry of the $j^{th}$ sample vector, and compute
the mean vector (in $\bfR^m$)
$$ \mbox{\boldmath${\mu}$} = \frac{1}{n}\sum_{k=1}^n \vect{d}_k,$$
\item Build the corresponding mean-centered data matrix $B$ with columns given by $\vect{d}_j - \mbox{\boldmath${\mu}$}$ so that the entries are
$$b_{ij} = x_{ij} - \mu_i$$ 
for every $i=1,...,m$ and $j=1,...,n$.
\item Use $B$ to compute the symmetric, $m \times m$ covariance matrix
$$ S = \frac{1}{n-1}B B^T.$$

\item Find the eigenvalues $\lambda_1,..., \lambda_m$ of $S$ (arranged in decreasing order including multiplicity) and an orthonormal set of \edits{corresponding} eigenvectors
$\vect{v}_1, ..., \vect{v}_m$. These create a new basis for $\bfR^m$ in which the data can be expressed. 

\item Finally, the data is represented in terms of the new basis vectors $\vect{v}_1, ..., \vect{v}_m$ using the coordinates $\vect{y}_1 = \vect{v}_1^T X, ..., \vect{y}_m = \vect{v}_m^T X$.  This can also be represented as the matrix $Y = P^TX$ where $P$ is the matrix with columns $\vect{v}_1, ..., \vect{v}_m$.  Should we wish to convert the data back to the original basis, we merely utilize the orthogonality of $P$ and compute $PY$ to find $X$.

\end{enumerate}

In the final sections, we will extend our introductory example while presenting additional applications in which PCA appears prominently \edits{and can} be implemented within a classroom environment. \edits{For each of the subsequent examples, an instructor could integrate the content into handouts or group worksheets, assign group projects with brief presentations during class, or merely use the given materials to present the information in an interactive lecture format.}

\section{Application - Data Analysis}
\label{data}

In the previous section, we developed a method for principal component analysis which determined \edits{a} basis with maximal variance.  Notice that the first component really encapsulated the majority of the information
embedded within the data.
Since the eigenvalues can be ordered, we might also
be able to truncate the sum in (\ref{ST}) to reduce the amount of stored data.
For instance, in our previous example, we might only keep the first principal component since
the data can be mostly explained just by knowing this characteristic, rather than every height and weight.
\edits{In this case, each data point would then be represented by its projection onto the first principal component.}
Upon performing the final step, we might also interpret the results: are a small number of the eigenvalues $\lambda_k$ \edits{much less (perhaps by an order of magnitude)} than the others? If so, this indicates a reduction in the dimension of the data is possible without losing much information, \edits{while}
if this does not occur then the dimension of the data may not be easily reduced \edits{in this way}.

Suppose that in addition to computing the components, we were to truncate
the \edits{new} basis matrix, $P$, so that we keep only the first $r$ columns with $r < m$. We would thus have a matrix
$\tilde{P} \in \bfR^{m\times r}$\edits{, and this would give rise to the $r \times n$ matrix $\tilde{Y} = \tilde{P}^T X$ containing
a truncation of the data represented by the first $r$ principal components.
From this we could also create the $m \times n$ matrix $\tilde{X} = \tilde{P}\tilde{P}^T X$, which represents the projection of $X$ onto the first $r$ principal directions.  This reduced representation of the data would then possess less information than $X$, but retain the most information when compared to any other matrix of the same rank.}
In fact, an error estimate is obtained from the Spectral Theorem - namely, the amount of information retained is given
by the {\em spectral ratio} \edits{of the associated covariance matrix}
\begin{equation}
\label{sigma}
\sigma^2  \newedits{:=}  \frac{ \sum_{k=1}^r \lambda_k}{ \sum_{k=1}^n \lambda_k}.
\end{equation}
Thus, in our first example, we can keep only the first component of each data point (a $1 \times 30$ matrix)
rather than the full data set ($2 \times 30$ matrix) and still retain $99\% $ of the information contained
within because
$$ \sigma^2 = \frac{191.46}{191.46 + 0.76} > 0.99.$$
In situations where the dimension of the input vector is large, but the components of the vectors are highly correlated, it is beneficial to reduce the dimension of the data matrix using PCA. This has three effects - it orthogonalizes the \edits{basis vectors} (so that they are uncorrelated), orders the resulting orthogonal components so that those with the largest \edits{variance appear} first, and eliminates dimensions that contribute the least to the variation in the data set.

\begin{figure}[t]
\begin{subfigure}[t]{\textwidth}
\centering
\includegraphics[scale=0.3]{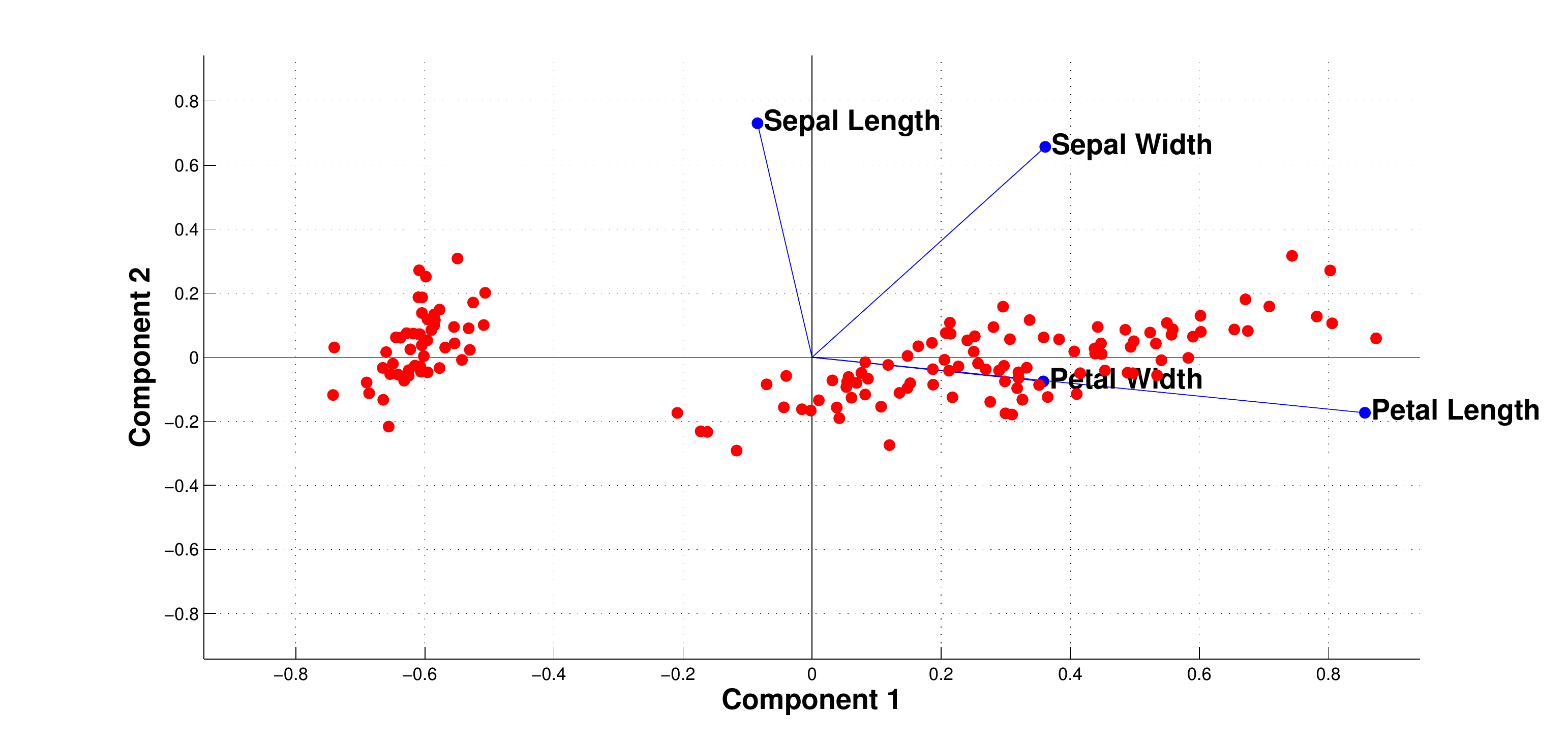} 
\vspace{-0.1in}
\caption{\footnotesize Biplot of Iris data.  }
\label{fig:BP2D}  
\end{subfigure}
\begin{subfigure}[t]{\textwidth}
\centering
\includegraphics[scale=0.3]{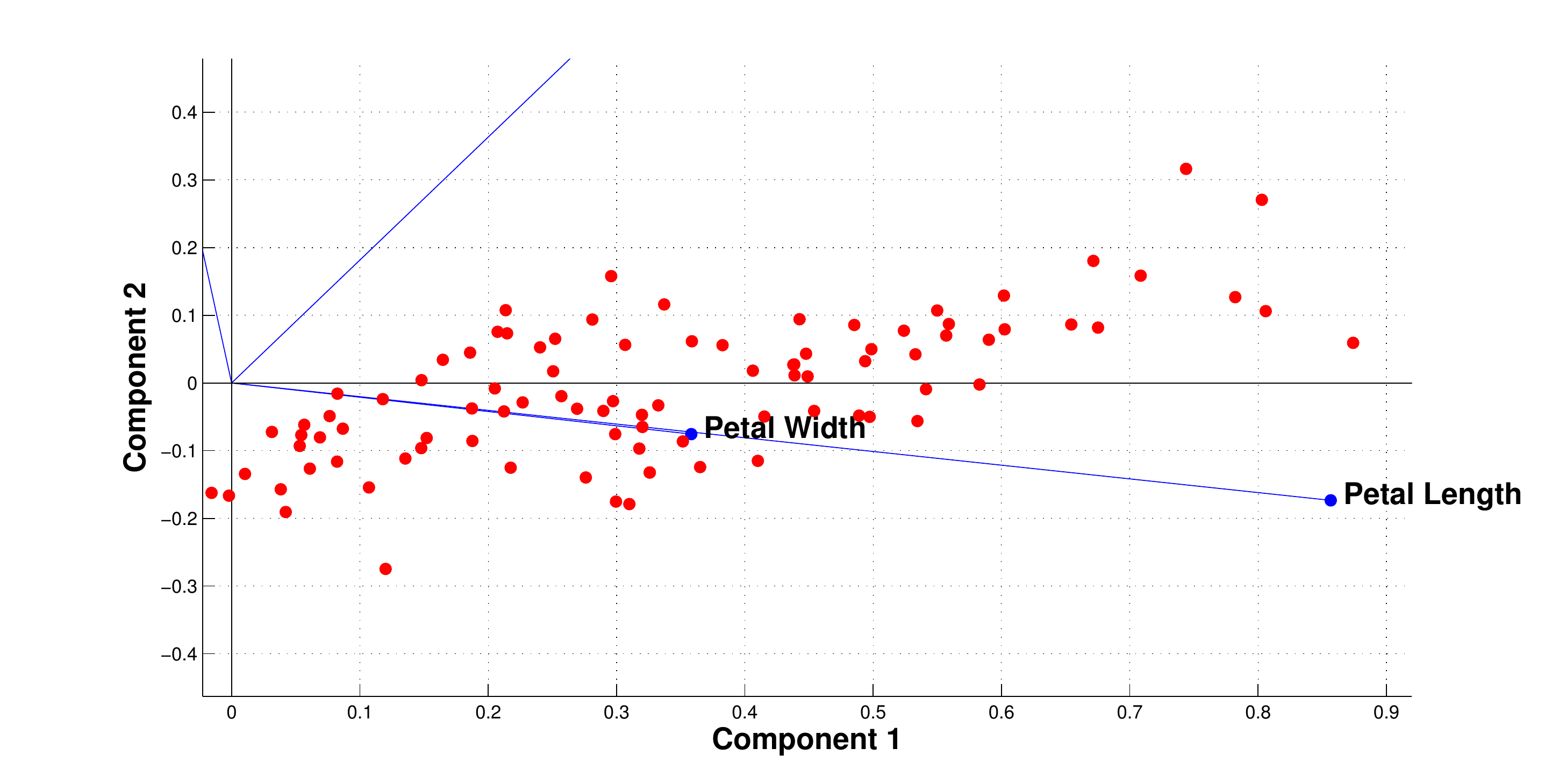} 
\vspace{-0.1in}
\caption{\footnotesize Enlarged portion of (a).  Notice that Petal Width and Petal Length point in nearly identical directions}
\label{fig:BP2DZoom} 
\end{subfigure}
\caption{Biplots of Fisher's Iris data projected onto $r = 2$ principal components with $\sigma^2 = 97.77\%$}
\label{Biplot}
\end{figure}

We now extend our introductory example to larger data sets with many other characteristics.
For our first example, we'll utilize the built-in {\it MATLAB} data set \texttt{fisheriris.mat}.
This famous collection of data arises from Fisher's 1936 paper \cite{Fisher} describing $50$ different samples of $4$ characteristics from each of three species of Iris.  Hence, the data set contains $150$ points and $4$ variables: \textbf{Sepal length}, \textbf{Sepal width}, \textbf{Petal length}, and \textbf{Petal width}.
With the $4 \times 150$ data matrix loaded, we perform the steps outlined within Section \ref{sec:steps} and compute the first two \edits{principal} components, i.e. those corresponding to the largest eigenvalues $\lambda_1 = 4.23$ and $\lambda_2 = 0.24$.
The others $\lambda_3 = 0.08$ and $\lambda_4 = 0.02$ are omitted.

It may be difficult to visualize even the first two principal components in this example because they are vectors in $\bfR^4$, but we can list them: 
\begin{table}[H]
\centering
\begin{tabular}{l | c c}
Characteristic & PC1 & PC2\\
\hline
Sepal Length  &	 0.3614   &  0.6566\\
Sepal Width 	&     	   -0.0845   & 0.7302\\
Petal Length	 &             0.8567    & -0.1734\\
Petal Width 	&            0.3583   & -0.0755
\end{tabular}
\label{tab:3}
\end{table}

\vspace{-0.2in}
\noindent Alternatively, we can visualize each data point projected onto the first two principal components as in  Fig.~\ref{fig:BP2D}.
This figure is actually a biplot, containing both the projected data points and 
the proportion of each characteristic which accounts for the respective principal component.  For instance, because the Petal length accounts for
a large proportion of the first principal component, it points in nearly the same direction as the $x$ axis within the biplot.  Similarly, Sepal width and Sepal length seem
to account for a large amount of the second principal component.

\begin{figure}[t]
\begin{subfigure}[t]{0.45\textwidth}
\hspace{-0.3in}
\includegraphics[scale=0.33]{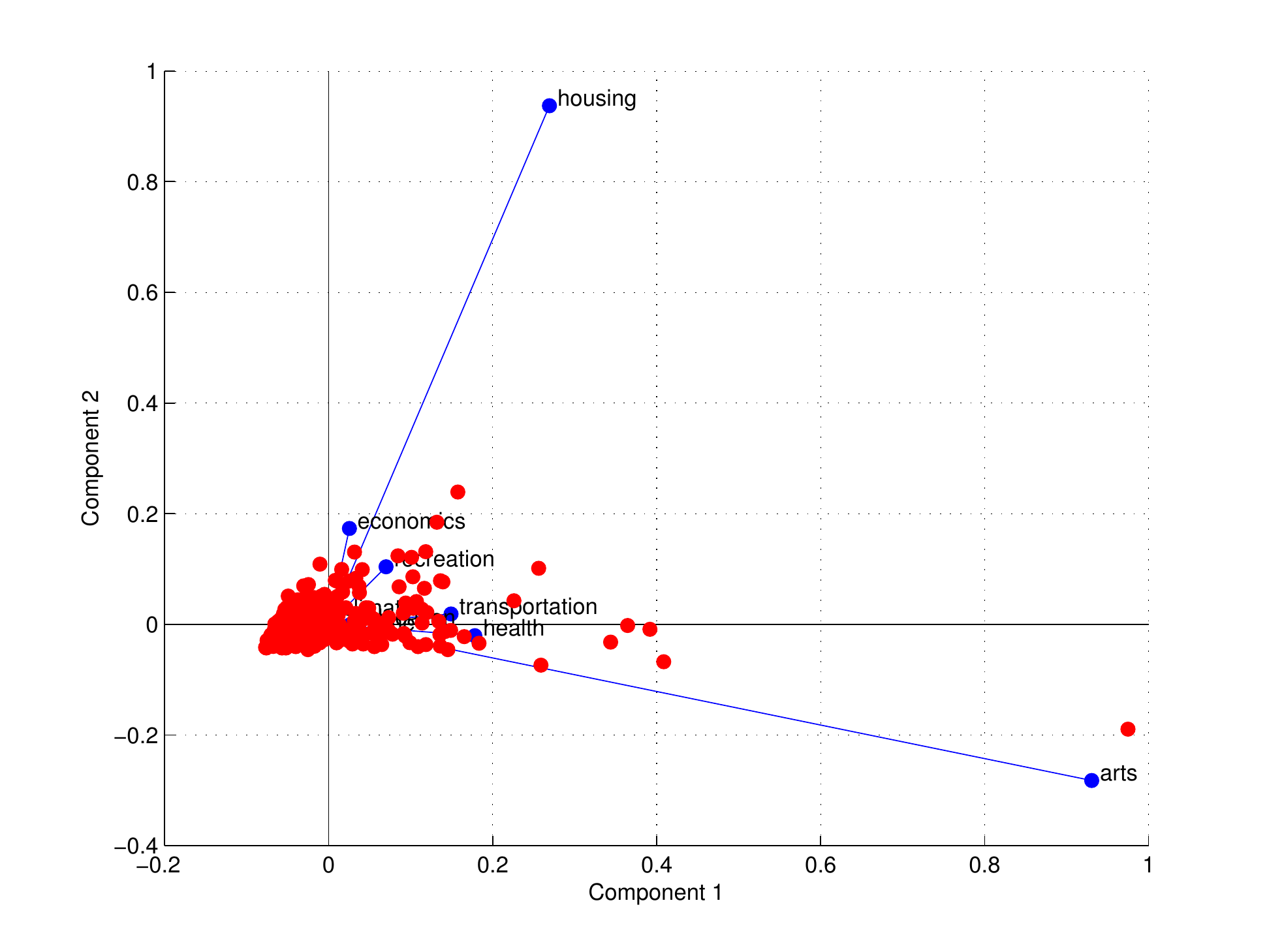}
\vspace{-0.3in}
\caption{\footnotesize $r=2$ components}
\label{fig:BPcities2D}  
\end{subfigure}
\hspace{0.1in}
\begin{subfigure}[t]{0.45\textwidth}
\includegraphics[scale=0.33]{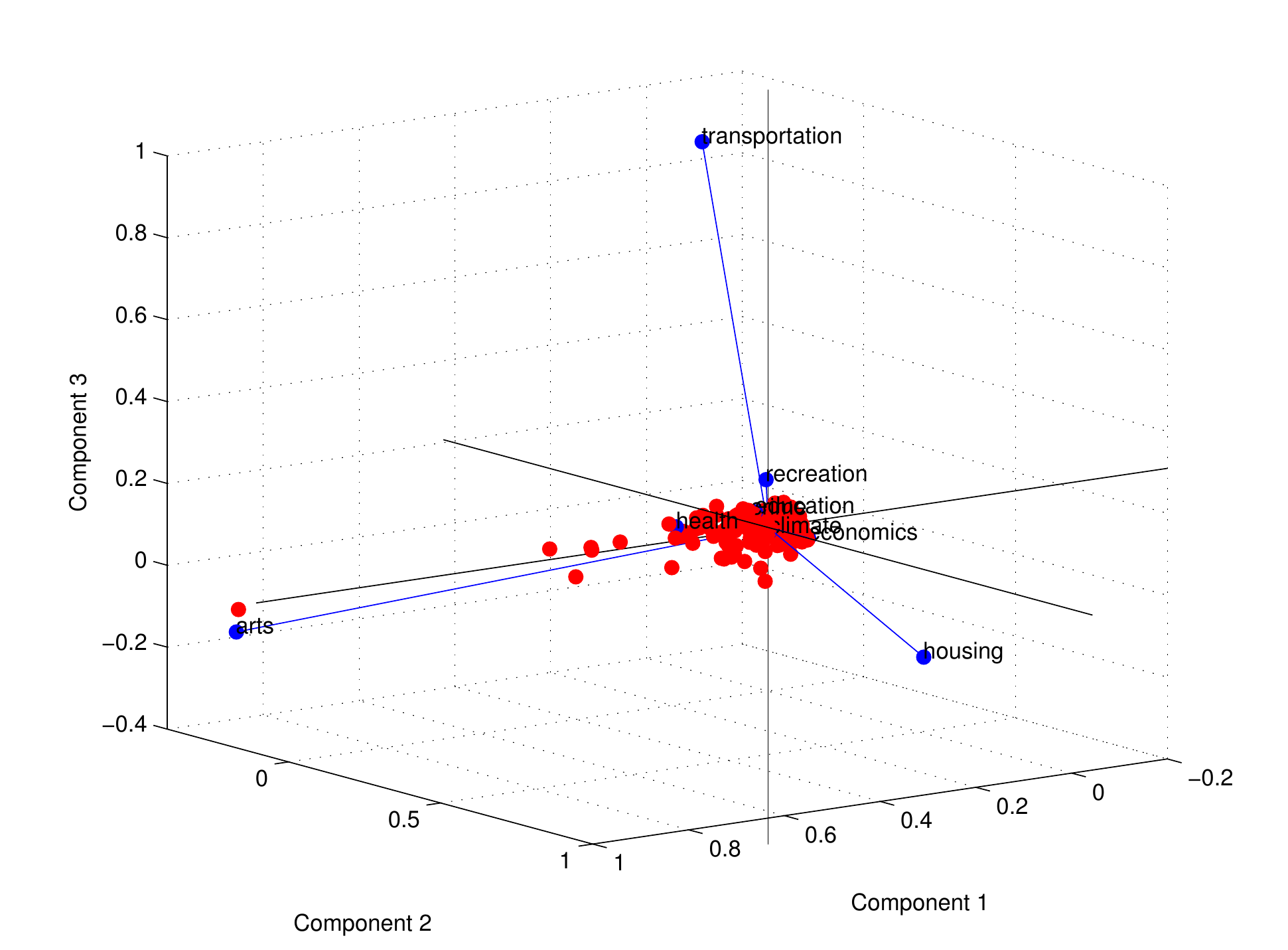}
\vspace{-0.3in}
\caption{\footnotesize $r=3$ components}
\label{fig:BPcities3D} 
\end{subfigure}
\caption{Biplots of \texttt{cities} data, both with $\sigma^2 \approx 94.43\%$}
\label{Biplot_cities}
\end{figure}

Notice in Fig.~\ref{fig:BP2D} that most of the data lie along the first component, which is mostly determined by the length of the Iris petals.  Additionally, Petal length and Petal width appear to be strongly correlated because they point in nearly identical directions within the new basis.  This can be seen in the enlarged portion of the biplot shown in  Fig.~\ref{fig:BP2DZoom}.  Contrastingly, Sepal width and Sepal length seem only mildly correlated and neither appears to be correlated to the pedal characteristics since \edits{these vectors point in somewhat orthogonal directions}.

Of course, the dimension and complexity of this example can be increased by using a larger set of data, such as another built-in {\it MATLAB} sample called \texttt{cities.mat}.  For completeness, we've included biplots (Fig.~\ref{Biplot_cities}) of the first few principal components of the \texttt{cities} set, which contains $m=9$ different attributes (i.e., climate, crime, education, etc..) for $n=329$ cities, again using the code contained in the appendix.
As a final observation, we note that due to the \edits{separation between eigenvalues of the covariance matrix (sometimes called the ``spectral gap'')}, the majority of data points lie along the first principal component in Fig.~\ref{Biplot_cities}(a) and within the plane generated by the first two principal components in Fig.~\ref{Biplot_cities}(b). \edits{This indicates that the majority of the variance within the data exists in these two directions, and hence one might safely eliminate the other dimensions, which contribute the least to the variation in the data set.}

\section{Application - Neuroscience}

\begin{figure}[t]
\centering
\includegraphics[scale=0.4]{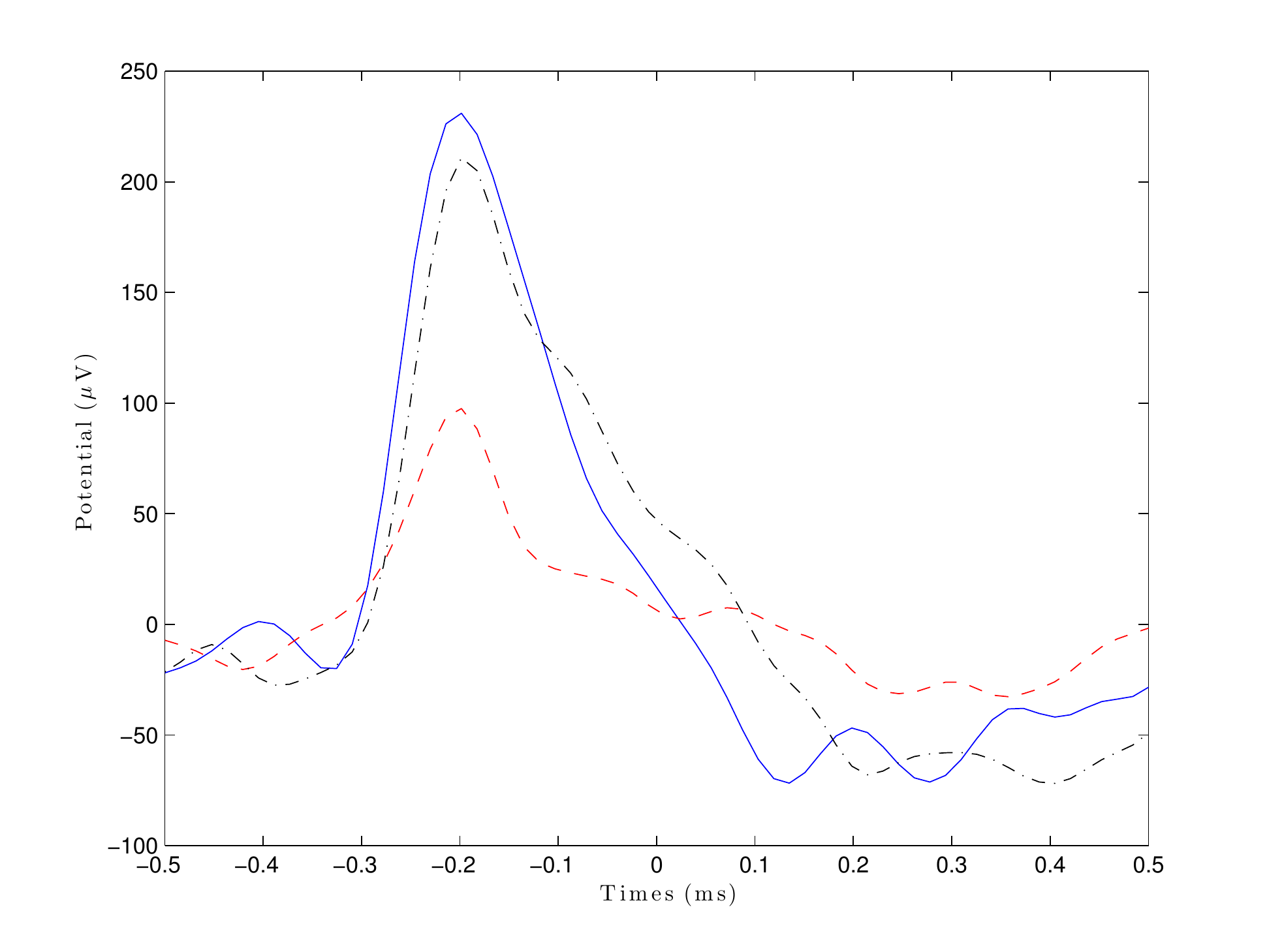}
\vspace{-0.2in}
\caption{ \footnotesize \label{fig:ap} Examples of action potentials}
\end{figure}

For another field in which PCA is quite useful, we turn to Neuroscience. In electrophysiological recordings of neural activity, one electrode typically records the action potentials (or spikes) of several neurons. 
Before one can use the recorded spikes to study the coding of information in the brain,
they must first be associated with the neuron(s) from which the signal arose.  
This is often accomplished by a procedure called \textit{spike sorting}, which can be accomplished because the recorded spikes of each neuron often have characteristic shapes \edits{\cite{spikes}}. 
For example, Fig.~\ref{fig:ap} shows three different shapes of action potentials recorded with an electrode.  
These potentials are plotted by connecting measurements at $64$ different points in time using linear interpolation.
\newedits{Due to distinctions in their peaks and oscillations}, the spikes in Fig.~\ref{fig:ap} are due to three different neurons, and PCA can be used to identify the principal variations of these spike shapes. 

\begin{figure}[t]
\hspace{-0.7in}
\includegraphics[scale=0.38]{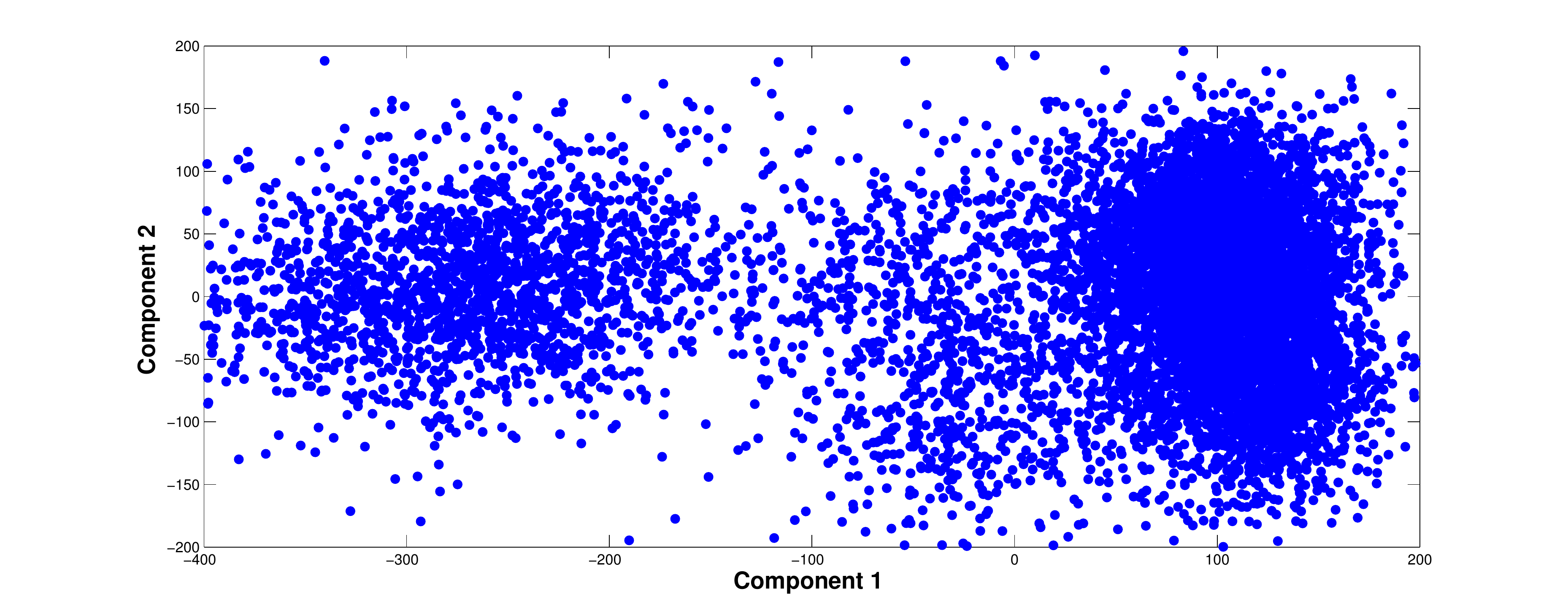}
\vspace{-0.4in}
\caption{ \footnotesize \label{fig:cluster} Spike data projected onto principal components}
\end{figure}

\begin{figure}[t]
\centering
\includegraphics[scale=0.36]{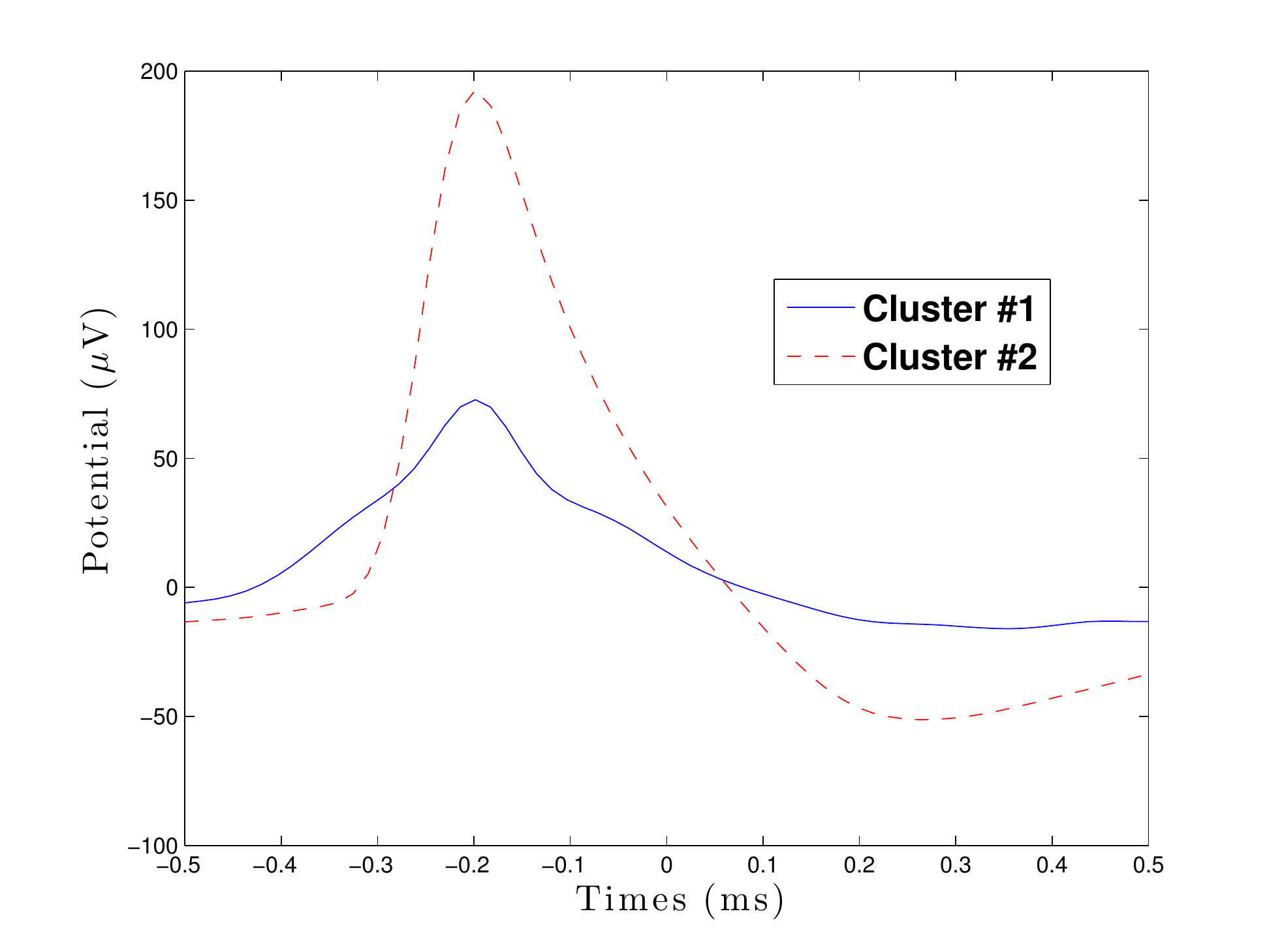}
\vspace{-0.2in}
\caption{ \footnotesize \label{fig:cluster_comp} Representative action potentials for clusters in projected data}
\end{figure}

Let's consider a data matrix $X$ of $9195$ recorded spikes, sampled at $64$ different points in time as in the previous figure
so that $X \in \mathbb{R}^{64 \times 9195}$.
For this example, we've used real human data that is freely-available online \cite{spikes}\newedits{, but note that this data is not directly connected to the examples of action potentials in Fig.~\ref{fig:ap}}. 
Using PCA we compute the eigenvalues of the covariance matrix to find that the first two 
($\lambda_1 = 2.9 \times 10^4$ and $\lambda_2 = 4.7\times 10^3$) account for $83\%$ of the total information.
\newedits{Of course, additional components can be included to increase this percentage.}  
A plot of the data projected along the first two principal components is given in Fig.~\ref{fig:cluster}.
We notice that two distinct clusters have formed within the data, and thus it appears that the spikes are formed from two
different neurons.  
If we denote the first two principal vectors by $\vect{v}_1$ and $\vect{v}_2$, then Cluster \#1 (right) and Cluster \#2 (left) appear to center around
$100\vect{v}_1 + 25\vect{v}_2$ and $-275\vect{v}_1$, respectively. The spike shapes corresponding to these vectors are displayed in 
Fig.~\ref{fig:cluster_comp} and represent averaged activity from the two neurons.  Hence, we see
that PCA can be used both to determine the degree of correlation amongst certain characteristics and to identify clustering patterns
within data.
Additionally, it provides us with a lower-dimensional picture of high-dimensional data sets, and this is often very useful when attempting
to visualize high-dimensional data.
With the action potentials determined, they can be associated to specific neurons and analyzed within neuroscience studies.

\section{Application - Image Compression}

Another important application of PCA is Image Compression.  Because images are stored as large matrices with real entries, one can reduce their storage requirements by keeping only the essential portions of the image \edits{\cite{Lay}}.  Of course, information (in this case, fine-grained detail of the image) is naturally lost in this process, but it is done so in an optimal manner, so as to maintain the most essential characteristics.

\begin{figure}[t]
\centering
\includegraphics[scale=0.4]{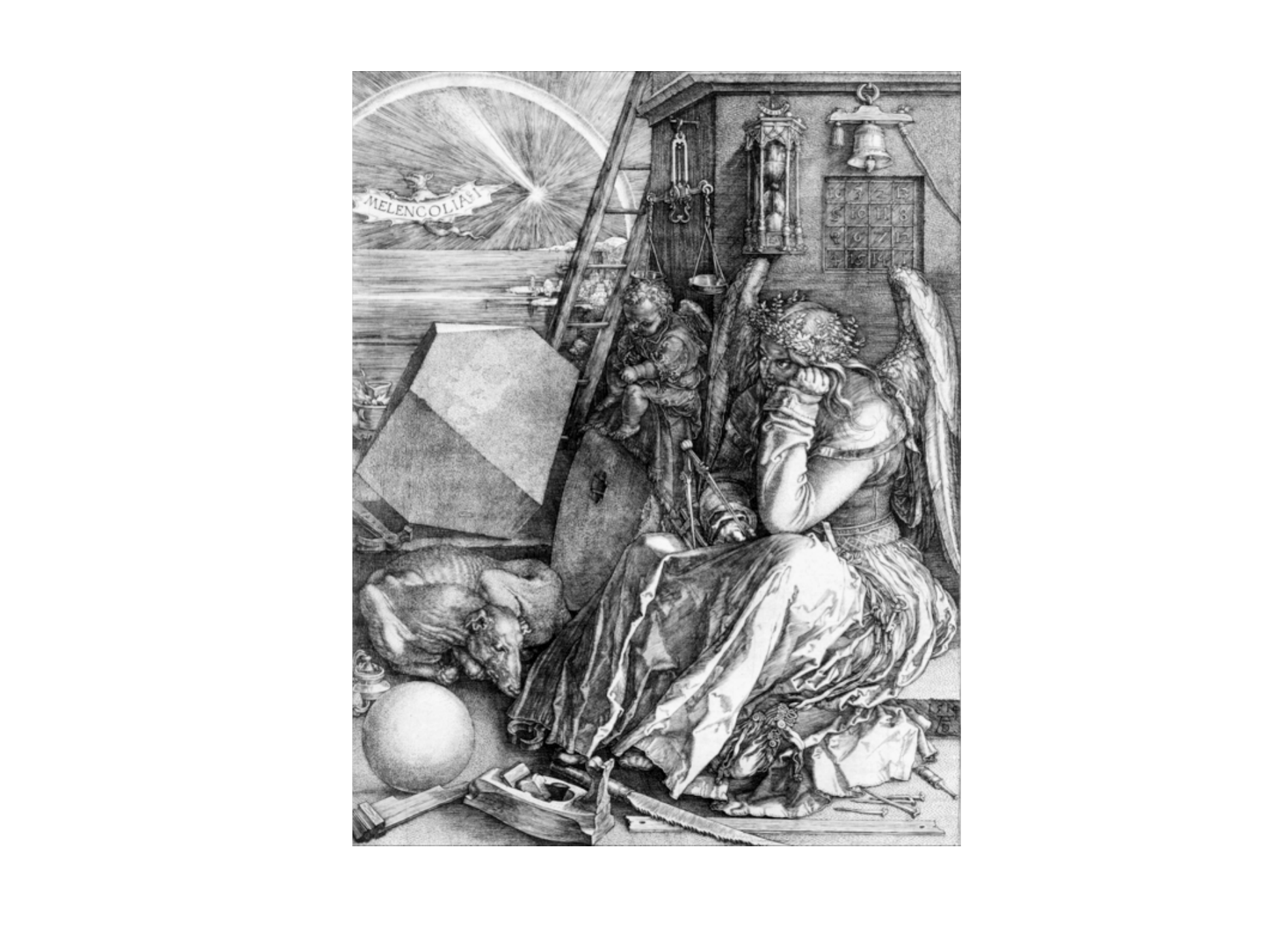}
\vspace{-0.2in}
\caption{ \footnotesize \label{fig:Durer} Albrecht D\"{u}rer's {\em Melancolia} displayed as a $648 \times 509$ pixelated imaged, taken from {\it MATLAB}'s built-in ``Durer'' file}
\end{figure}

\begin{figure}[t]
\centering
\includegraphics[scale=0.35]{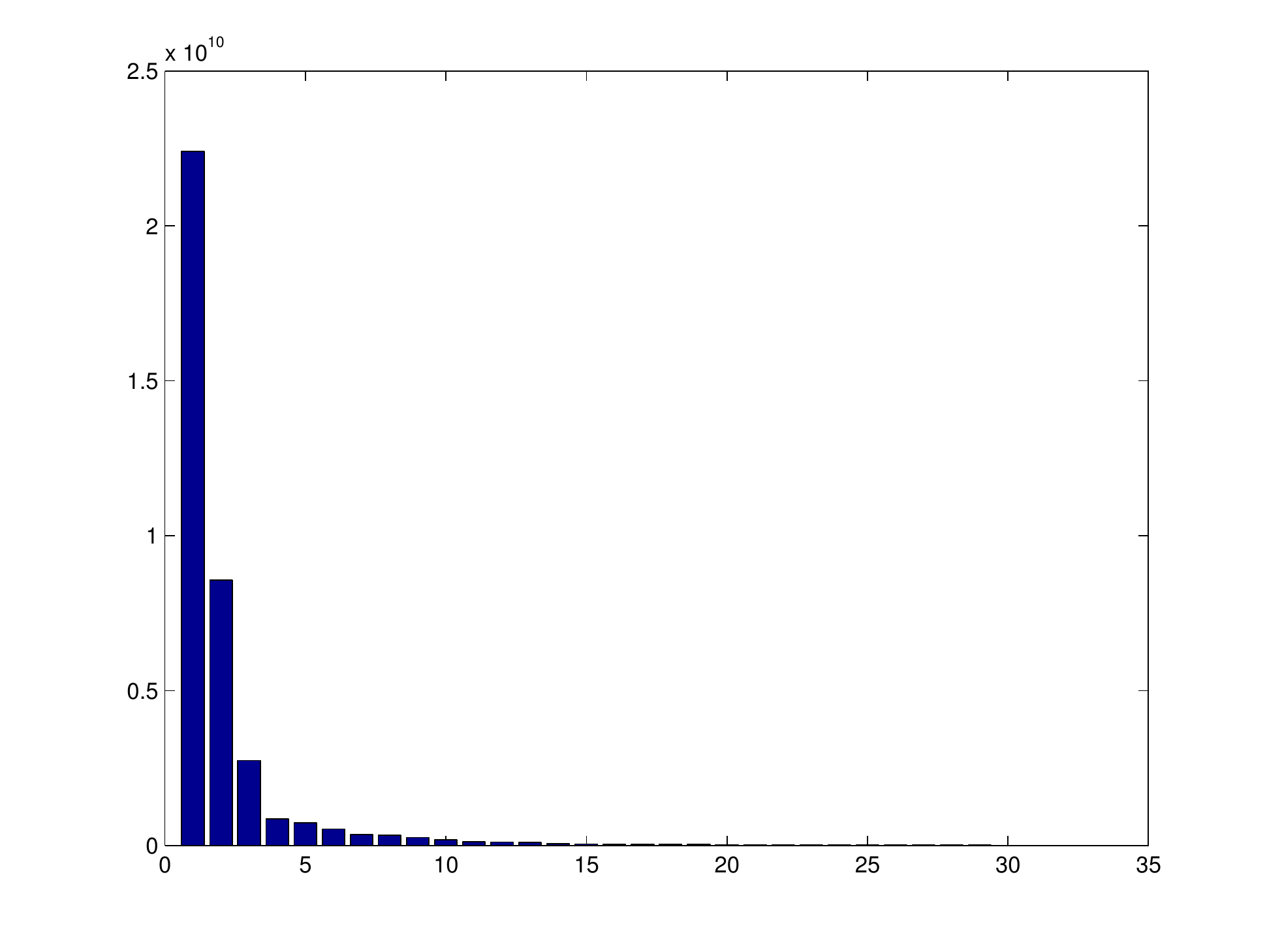}
\vspace{-0.2in}
\caption{ \label{fig:eigen}  \footnotesize  The first $35$ eigenvalues of \edits{the covariance matrix $S$ generated from $X$ in the image compression example and} arranged in decreasing order}
\end{figure}

\begin{figure}[t]
\begin{subfigure}[t]{0.45\textwidth}
\hspace{-0.3in}
\includegraphics[scale=0.3]{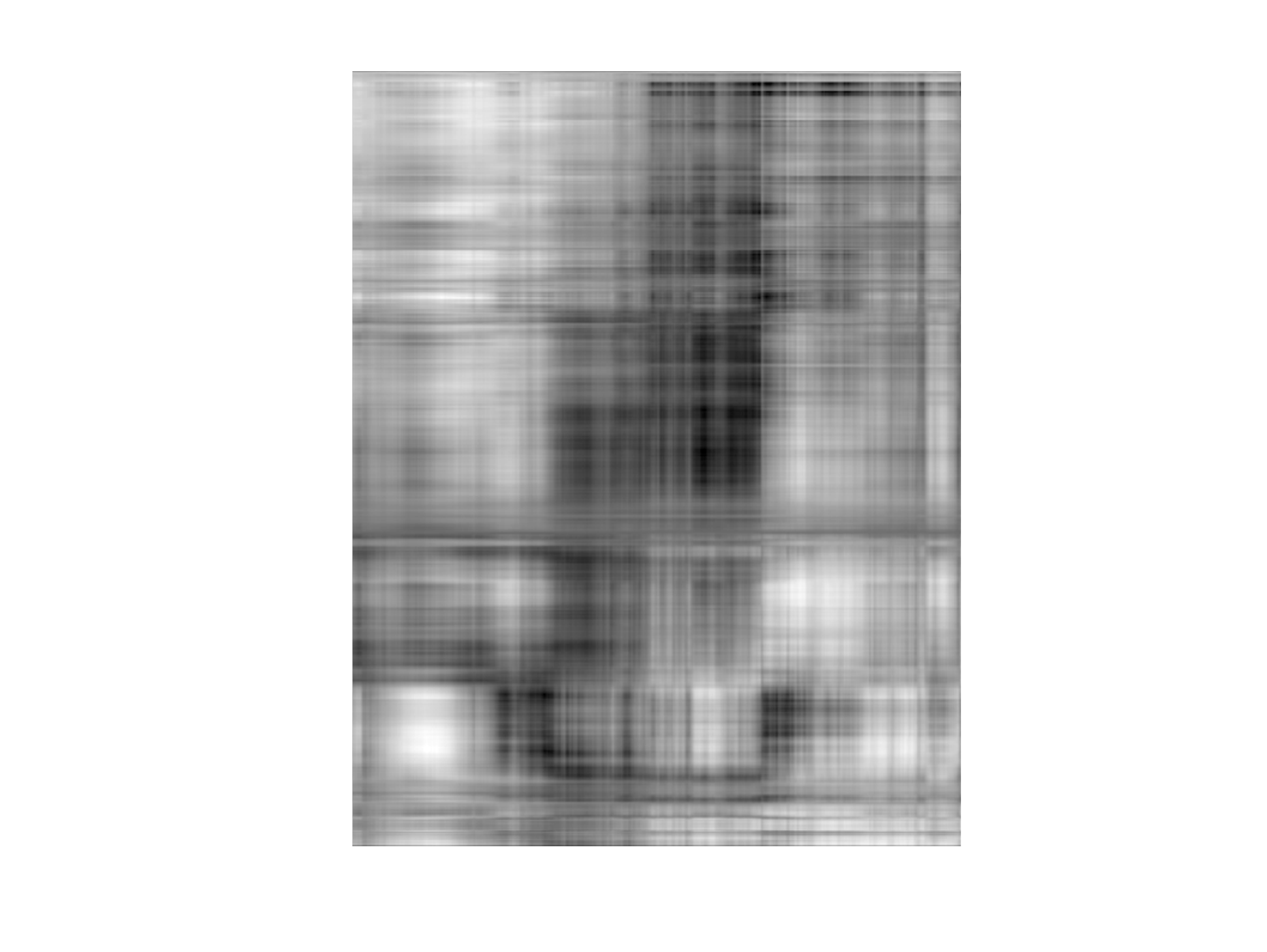}
\vspace{-0.3in}
\caption{ \label{figure:PC3}  \footnotesize  $r = 3$, $\sigma^2 = 88.89\%$}
\end{subfigure}
\hspace{0.1in}
\begin{subfigure}[t]{0.45\textwidth}
\hspace{-0.2in}
\includegraphics[scale=0.3]{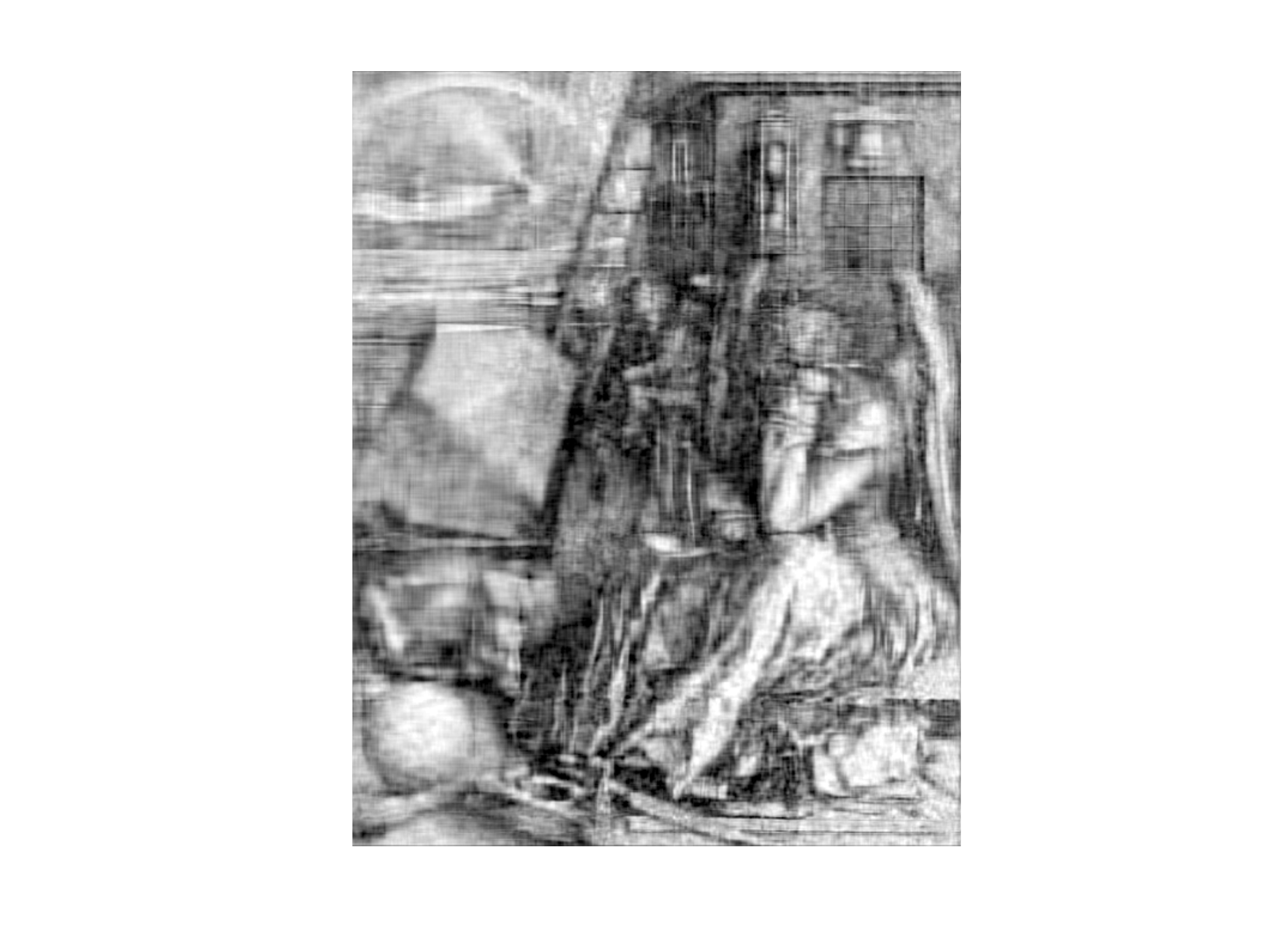}
\vspace{-0.3in}
\caption{ \label{figure:PC30}  \footnotesize $r = 30$, $\sigma^2 = 99.61\%$}
\end{subfigure}
\\
\begin{subfigure}[t]{0.45\textwidth}
\hspace{-0.3in}
\includegraphics[scale=0.3]{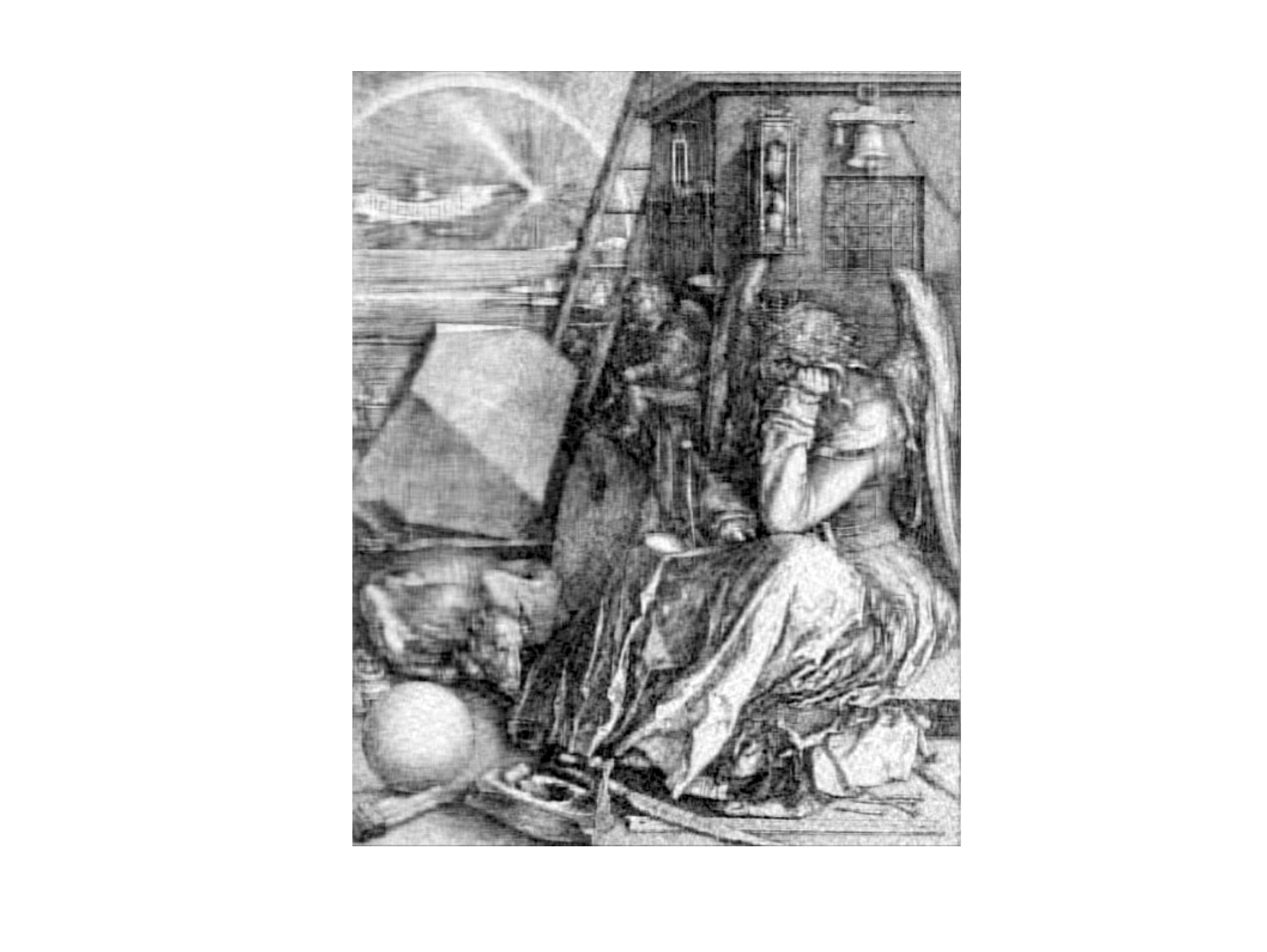}
\vspace{-0.3in}
\caption{ \label{figure:PC60}  \footnotesize $r = 60$, $\sigma^2 = 99.78\%$}
\end{subfigure}
\hspace{0.1in}
\begin{subfigure}[t]{0.45\textwidth}
\hspace{-0.2in}
\includegraphics[scale=0.3]{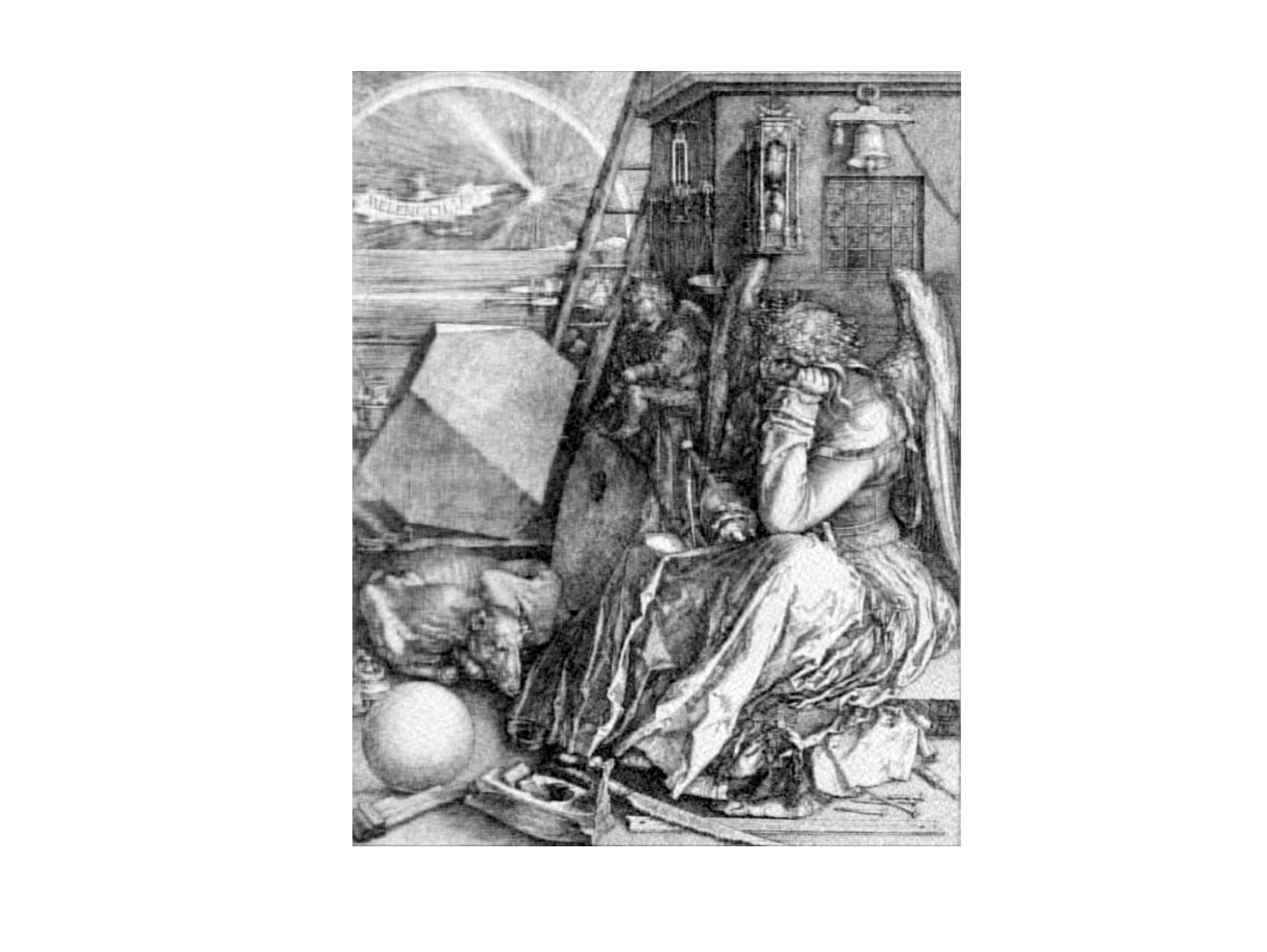}
\vspace{-0.3in}
\caption{ \label{figure:PC90}  \footnotesize  $r = 90$, $\sigma^2 = 99.92\%$}
\end{subfigure}
\caption{The D\"{u}rer image with varying \edits{numbers of} principal components.}
\label{Durer}
\end{figure}

\edits{In this section we detail a specific example for the use of PCA to compress an image.  
Since the effects of keeping a lower dimensional projection of the image will be visually clear, this particular example
is a great candidate for an interactive, in-class activity.  
More specifically, one can provide students with the code given in the appendix and ask them to determine the number of principal components $r$ that they must preserve in order to visually identify the image.  Additionally, since the variance in the projection of the original image onto the reduced basis is computed in the code, students could be asked to identify the value of $r$ that is needed to capture a certain percentage of the total variance. For instance, Fig. \ref{Durer} shows that $r$ must be at least $60$ in order to capture $99.78\%$ of the image detail.} 

\edits{Throughout the example} we will work with a built-in test image - 
Albrecht D\"{u}rer's {\em Melancolia} displayed in Fig.~\ref{fig:Durer}. 
{\it MATLAB} considers greyscale images like this as objects consisting of two portions - a matrix of
pixels and a colormap. 
Our image is stored in a $648 \times 509$ pixel matrix, and
thus contains $648 \times 509 = 329,832$ total pixels. 
The colormap is a $648 \times 3$ matrix, which we will ignore for the current study. 
Each element of the pixel matrix contains a real number representing the intensity of grey
scale for the corresponding pixel. 
{\it MATLAB} displays all of the pixels simultaneously
with the correct intensity, and the greyscale image that we see is produced.
The $648 \times 509$ matrix containing the pixel information is our data matrix, $X$.
Because the most important information in the reduced matrix $\tilde{X}$, described in Section \ref{data}, is captured by the first few principal components this suggests a 
way to compress the image by \edits{using the lower-rank approximation} $\tilde{X}$.


Computing the distribution of associated eigenvalues as in previous examples, we see the formation of a large spectral gap, as shown in Fig.~\ref{fig:eigen}.  Hence, the truncation of the sum of principal components in $\tilde{X}$ should still contain a large amount of the total information of the original image $X$.
In Fig.~\ref{Durer}, we've represented $\tilde{X}$ for four choices of $r$ (i.e., the number of principal components \edits{used}), and the associated
spectral ratio, $\sigma^2$, retained by those reduced descriptions 
is also listed.  
Notice that the detail of the image improves as $r$ is increased,
and that a fairly suitable representation can be obtained with around $90$ components rather than the full $648$ vector description.
Hence, PCA has again served the useful purpose of reducing the dimension of the original \newedits{data set} while preserving its most essential features.

%
%
\section{CONCLUSION}
The aim in writing this article is to generate and present clear, concise resources for the 
instruction of Principle Component Analysis (PCA) within \edits{an introductory} Linear Algebra course.
Students should be led to understand that PCA is a powerful, useful tool that is utilized \edits{in data analysis and} throughout
the sciences.
The visual components of the Image Compression examples in particular are designed to display the utility of PCA
and linear algebra, in general.
\edits{Within a more advanced setting (perhaps a second semester course in linear algebra),} new material regarding the Singular Value Decomposition can also be included to \edits{introduce additional details concerning the implementation of} PCA \edits{and provide an elegant structure for performing the method} \cite{Kalman, TB}. That being said, the SVD can often be a technical and 
time-consuming topic that one first encounters in a computational\edits{, honors, or graduate} linear algebra course
rather than within an introductory setting.
In general, we believe that these instructional resources should be helpful in seamlessly integrating one of the most essential current
applications of linear algebra, PCA, into a standard undergraduate course.  

\section*{ACKNOWLEDGEMENTS}
The authors would like to thank the College of Engineering and Computational Sciences at the Colorado School of Mines for partial support to develop and implement the educational tools included within this article.
Additionally, the first author thanks the National Science Foundation for  support under award DMS-1211667.

\newpage

\section*{APPENDIX - Matlab Code for Examples}

%
%

\vspace{-0.1in}

\noindent Load Height/Weight, Iris, cities, Neuroscience, or Durer image data and call functions. \newedits{Recall that the 
``HtWtdata'' and ``times\_CSC4'' data sets cannot be loaded directly from MATLAB, but can be obtained online using the aforementioned references.}
\vspace{-0.2in}
\begin{lstlisting}
%Initialize data for one of the five examples below

%%% For Height/Weight example, use
% load HtWtdata 	
% vbls = {'Height','Weight'};
% vbls=[];
% plottype = 1;
  
%%% For Iris data example, use 
% load fisheriris
% data = meas';
% vbls = {'Sepal Width','Sepal Height','Petal Length','Petal Width'};
% plottype = 2;

%%% For Cities example, use
% load cities; data = ratings'; vbls = categories;
% plottype = 2;

%%% For Durer example, use
% load durer; data=double(X)
% fig1=figure('Color',[1 1 1]);
% image(X), colormap(map), axis off, axis equal; % plot
% vbls=[]; plottype = 3;

%%% For Neuroscience example, use
% load times_CSC4
% data = spikes';
% vbls=[];
% plottype = 4; % plottype = 5;

PC = 2;	%Specify number of components to use
%Function calls for PCA and plotting
[X, V, Y] = PCAfun(data, PC);
PlotData(plottype, X, V, Y, PC, vbls, data);
\end{lstlisting}

\newpage

\noindent Function to plot components, biplots, and images
\vspace{-0.2in}
\begin{lstlisting}
function PlotData(num, X, V, Y, PC, vbls, data)

fig3 = figure('Color',[1 1 1]); 
if num == 1   		 % Plot for Statistics example
 plot(X(1,:), X(2,:),'.', 'MarkerSize', 30)
 xlabel('Height (in)')
 ylabel('Weight (lbs)')
 
elseif num == 2 	% Biplot for Statistics example
    Y = Y';
    biplot(V,'scores',Y(:,1:PC), 'MarkerSize', 20,'varlabels',vbls)

elseif num == 3		% Plot for Image Compression example
    image(X),colormap(map),axis off, axis equal; % display results

elseif num == 4		% Plot for Neuroscience Cluster
 plot(Y(1,:), Y(2,:),'.', 'MarkerSize', 20), axis([-400, 200, -200, 200])
 xlabel('Component 1 ', 'FontSize', 10)
 ylabel('Component 2 ', 'FontSize', 10)
     
elseif num == 5		% Plot for Action Potentials
 time = linspace(-0.5, 0.5, 64);
 mu = mean(data,2);
 pot1 = mu + 110*V(:,1) - 25*V(:,2);
 pot2 = mu - 260*V(:,1) - 10*V(:,2);
 plot(time, pot1 , time,  pot2, '--r')
 xlabel('Times (ms)','Interpreter','latex')
 ylabel('Potential ($\mu$V)','Interpreter','latex')
 legend('Cluster #1', 'Cluster #2')
  
end

figureHandle = gcf;
%%% make all text in the figure to size 16 and bold
set(findall(figureHandle,'type','text'),'fontSize',
16,'fontWeight','bold')
\end{lstlisting}

\newpage

\noindent Function to compute principal components from given data
\vspace{-0.2in}
\begin{lstlisting}
function [Xtrunc, Vtrunc, Y] = PCAfun (data, PC)
% inputs load data and establish number of components to use

 [m n] = size(data); 	% obtain dimensions of data matrix
 mu = mean(data,2); 	% compute row mean
 X = data-repmat(mu,1,n); 	% subtract row mean to obtain X
 Z = 1/sqrt(n-1)*X'; 	% create matrix, Z
 covar = Z'*Z; 	% covariance matrix of Z

% Extract eigenvalues and eigenvectors
 [V, Lambda] = eigs(covar); % compute eigenvectors and eigenvalues
 evals = diag(Lambda); % store eigenvalues in vector form
 fig2 = figure('Color',[1 1 1]); 
 bar(evals);	% bar plot of eigenvalues (variances)
 
 % compute percentage of variance included
 totvals = sum(evals);	%total sum of eigenvalues
 runtot = cumsum(evals);	%cumulative sum of eigenvalues - largest to smallest
 percentinc = runtot(PC)/totvals 	% percentage included

% Extract principal components
 Vtrunc = V(:,1:PC); % extract first PC vectors in new basis
 Y = Vtrunc' * X; % project data onto new basis
 Xtrunc = Vtrunc*Y; 	% convert back to original basis
 Xtrunc = Xtrunc+repmat(mu,1,n); 	% add the row means back
 \end{lstlisting}

\newpage

\section*{BIOGRAPHICAL SKETCHES}

Prof. Stephen Pankavich received his Ph.D. in Mathematical Sciences from Carnegie Mellon University in $2005$ and has been interested in applications of mathematics throughout his career. He is currently an Assistant Professor at the Colorado School of Mines where he actively seeks to relate mathematics to the daily life of his undergraduate students while stressing it's immense utility within science and engineering.  

\vspace*{.3 true cm} 

\noindent Prof. Rebecca Swanson received her Ph.D. in Mathematics from Indiana University in $2010$. Prior to arriving at the Colorado School of Mines in $2012$, she was an Assistant Professor at Nebraska Wesleyan University.  She enjoys the beauty and elegance of discrete mathematics and the challenge of pedagogical innovation.  In her free time, she also enjoys baking, running, and reading.


\end{document}